\documentclass[12pt]{article}
\usepackage{amsmath,amssymb,mathrsfs}
\usepackage[PostScript=dvips,s=5ex]{diag}
\usepackage{hyperref}
\usepackage[latin1]{inputenc}
\diagramstyle[heads=curlyvee,s=5ex]
\begin{document}
\title{A General Local-to-Global Principle for Convexity of Momentum Maps}

\vspace{-21mm}
\author{Wolfgang Rump and Jenny Santoso\\ {\footnotesize\it Institute for
Algebra and Number Theory, University of Stuttgart}\\[-2mm]
{\footnotesize\it Pfaffenwaldring 57, D-70550 Stuttgart, Germany}\\[4mm]
%
{\small\it Dedicated to Tudor S. Ratiu and Alan D. Weinstein}\\
{\small\it for their tireless \& outstanding support}\\ 
%
%
%
}

\date{Preprint 2009}
\maketitle
\setcounter{footnote}{-1}
\renewcommand{\thefootnote}{\Roman{footnote}}
\footnotetext{\hspace*{-1.5mm} 2000 \it Mathematics Subject Classification.
\rm Primary: 52A01, 53C23, 54C10, 53D20.
Secondary: \rm 54E18.

\it \hspace*{1mm} Key words and phrases. \rm Convexity space, momentum map,
geodesic manifold, Lokal-global-Prinzip.}
\thispagestyle{empty}

\newtheorem{prop}{Proposition}
\newtheorem{thm}{Theorem}
\newtheorem{lem}{Lemma}
\newtheorem{Definition}{Definition}
\renewcommand{\labelenumi}{\rm (\alph{enumi})}
\newcommand{\hs}{\hspace{2mm}}
\newcommand{\vsp}{\vspace{4ex}}
\newcommand{\vspc}{\vspace{-1ex}}
\newcommand{\hra}{\hookrightarrow}
\newcommand{\tra}{\twoheadrightarrow}
\newcommand{\md}{\mbox{-\bf mod}}
\newcommand{\Mod}{\mbox{-\bf Mod}}
\newcommand{\Mdd}{\mbox{\bf Mod}}
\newcommand{\mdd}{\mbox{\bf mod}}
\newcommand{\latt}{\mbox{-\bf lat}}
\newcommand{\Proj}{\mbox{\bf Proj}}
\newcommand{\Inj}{\mbox{\bf Inj}}
\newcommand{\Ab}{\mbox{\bf Ab}}
\newcommand{\CM}{\mbox{-\bf CM}}
\newcommand{\Prj}{\mbox{-\bf Proj}}
\newcommand{\prj}{\mbox{-\bf proj}}
\newcommand{\ra}{\rightarrow}
\newcommand{\eps}{\varepsilon}
\renewcommand{\epsilon}{\varepsilon}
\renewcommand{\phi}{\varphi}
\renewcommand{\hom}{\mbox{Hom}}
\newcommand{\ex}{\mbox{Ext}}
\newcommand{\rad}{\mbox{Rad}}
\renewcommand{\Im}{\mbox{Im}}
\newcommand{\oti}{\otimes}
\newcommand{\sig}{\sigma}
\newcommand{\en}{\mbox{End}}
\newcommand{\Qq}{\mbox{Q}}
\newcommand{\lra}{\longrightarrow}
\newcommand{\lras}{\mbox{ $\longrightarrow\raisebox{1mm}{\hspace{-6.5mm}$\sim$}\hspace{3mm}$}}
\newcommand{\Eq}{\Leftrightarrow}
\newcommand{\Equ}{\Longleftrightarrow}
\newcommand{\Ra}{\Rightarrow}
\newcommand{\A}{\mathbb{A}}
\newcommand{\N}{\mathbb{N}}
\newcommand{\Z}{\mathbb{Z}}
\newcommand{\R}{\mathbb{R}}
\newcommand{\K}{\mathbb{K}}
\newcommand{\F}{\mathbb{F}}
\newcommand{\B}{\mathfrak{B}}
\newcommand{\rk}{\mbox{r}}
\newcommand{\p}{\mathfrak{p}}
\newcommand{\q}{\mathfrak{q}}
\renewcommand{\k}{\mathfrak{k}}
\newcommand{\AAA}{\mathfrak{A}}
\newcommand{\Aa}{{\cal A}}
\newcommand{\Dd}{{\cal D}}
\newcommand{\Bb}{{\cal B}}
\newcommand{\Ll}{{\cal L}}
\newcommand{\Hh}{{\cal H}}
\newcommand{\Nn}{{\cal N}}
\newcommand{\Cc}{{\cal C}}
\newcommand{\Mm}{{\cal M}}
\newcommand{\Tt}{{\cal T}}
\newcommand{\Ff}{{\cal F}}
\newcommand{\Rr}{{\cal R}}
\renewcommand{\AA}{\mathscr{A}}
\newcommand{\DD}{\mathscr{D}}
\newcommand{\BB}{\mathscr{B}}
\newcommand{\LL}{\mathscr{L}}
\newcommand{\HH}{\mathscr{H}}
\newcommand{\NN}{\mathscr{N}}
\newcommand{\EE}{\mathscr{E}}
\newcommand{\CC}{\mathscr{C}}
\newcommand{\UU}{\mathscr{U}}
\newcommand{\MM}{\mathscr{M}}
\newcommand{\PP}{\mathscr{P}}
\newcommand{\II}{\mathscr{I}}
\newcommand{\TT}{\mathscr{T}}
\newcommand{\XX}{\mathscr{X}}
\newcommand{\YY}{\mathscr{Y}}
\newcommand{\KK}{\mathscr{K}}
\newcommand{\FF}{\mathscr{F}}
\newcommand{\RR}{\mathscr{R}}
\newcommand{\PPP}{\mathfrak{P}}
\newcommand{\BBB}{\mathfrak{B}}
\newcommand{\DDD}{\mathfrak{D}}
\newcommand{\UUU}{\mathfrak{U}}
\newcommand{\FFF}{\mathfrak{F}}
\newcommand{\VVV}{\mathfrak{V}}
\newcommand{\WWW}{\mathfrak{W}}
\newcommand{\setm}{\smallsetminus}
\renewcommand{\le}{\leqslant}
\renewcommand{\ge}{\geqslant}
\newcommand{\rat}{\rightarrowtail}
\newcommand{\op}{^{\mbox{\scriptsize op}}}
\newcommand{\pf}{\it Proof.\hs\rm}
\newcommand{\bx}{\hspace*{\fill} $\square$}
\newcommand{\ltra}{\lra\!\!\!\!\!\!\:\ra}
\newcommand{\IP}{\tbinom{I}{P}}
\newcommand{\lap}{\Lambda\!^{\raisebox{.5mm}{\tiny $\! +$}}}
\newcommand{\lam}{\Lambda\!^{\raisebox{.5mm}{\tiny $\! -$}}}
\newcommand{\upll}{^{\raisebox{.5mm}{\tiny $++$}}}
\newcommand{\dpl}{\raisebox{-.5mm}[-1ex][-1ex]{\tiny $+$}}
\newcommand{\dmi}{\raisebox{-.5mm}[-1ex][-1ex]{\tiny $-$}}
\newcommand{\subs}{\subset}
\newcommand{\sups}{\supset}
\newcommand{\subsn}{\subsetneq}
\newcommand{\nsubs}{\not\subset}
\newcommand{\pt}{\makebox[0pt][r]{\bf .\hspace{.5mm}}}
\newcommand{\KS}{Krull-Schmidt }
\renewcommand{\theta}{\vartheta}
\newcommand{\rrr}{\rule[1mm]{5mm}{.5mm}}
\newcommand{\rr}{\rule[1mm]{2.5mm}{.5mm}}
\newcommand{\smr}{\rule[.5mm]{1.5mm}{.3mm}}
\newcommand{\dpt}{\makebox[0mm][l]{.}}
\newcommand{\dpc}{\makebox[0mm][l]{,}}
\newcommand{\MA}{\mbox{\textsf{M}$(\AA)$}}
\newcommand{\radp}{\underline{\rad}}
\newcommand{\radi}{\overline{\rad}}

\newarrow{To} ----{->}
\newarrow{to} ----{>}
\newarrow{in} {>}---{>}
\newarrow{inc} C---{>}
\newarrow{Inc} C---{->}
\newarrow{on} ----{>>}
\newarrow{On} ----{->>}
\newarrow{eq} =====
\newarrow{dash} {}{dash}{}{dash}{>}
\newarrow{dsh} {}{dash}{}{dash}{}

{\small\bf Abstract.} {\small We extend the
Local-to-Global-Principle used in the proof of convexity theorems
for momentum maps to not necessarily closed maps whose target space
carries a convexity structure which need not be based on a metric.
Using a new factorization of the momentum map, convexity of its
image is proved without local fiber connectedness, and for almost
arbitrary spaces of definition.}
%

%

\vspace{5mm}
\noindent {\Large\bf Introduction}

\vspace{5mm}
Convexity for momentum maps was discovered independently by Atiyah \cite{At}
and Guillemin-Sternberg \cite{GS} in the case of a Hamiltonian torus action on
a compact symplectic manifold $X$. It was proved that the image of the momentum
map
$\mu$ is a convex polytope, namely, the convex hull of $\mu(X^T)$, where $X^T$
denotes the set of fixed points under the action of the torus $T$. In this
case, $\mu$ is
open onto its image, and the fibers of $\mu$ are compact and connected. Two
years later, in 1984, Kirwan \cite{Kir} (see also \cite{GS1}) extended this
result to the action of a compact connected Lie group $G$. Here the image of
$\mu\colon X\ra\mbox{Lie}(G)^\ast$ has to be restricted to a closed Weyl
chamber in a Cartan subalgebra of $\mbox{Lie}(G)$, i.~e. a fundamental domain
of $G$ with respect to its coadjoint action on $\mbox{Lie}(G)^\ast$.
Equivalently, this amounts to a composition of the momentum map $\mu$ with the
projection onto the quotient space $Y:=\mbox{Lie}(G)^\ast/G$ modulo the
coadjoint action of $G$.
Up to this time, convexity of $\mu$ was proved by means of Morse theory,
applied to the components of $\mu$. This works well as long as $\mu$ is
defined on a compact manifold~$X$.

\textheight22.9cm 

In 1988, Condevaux, Dazord, and Molino
\cite{CDM} reproved these results in an entirely new fashion. They factor out
the connected components of the fibers of $\mu$ to get a monotone-light
factorization $\mu\colon X\ra\widetilde{X}\ra Y$ (see \cite{Mic}). If $\mu$ is
proper, i.~e. closed and with quasi-compact fibers, the
metric of $Y$ can be lifted to $\widetilde{X}$. Then a shortest path between
two points of $\widetilde{X}$ maps to a straight line in $Y$,
which proves the convexity of $\mu(X)$. Based on this method,
Hilgert, Neeb, and Plank \cite{HNP} extended Kirwan's result to non-compact
connected manifolds $X$ under the assumption that $\mu$ is proper.

After this invention, the proof of convexity now splits into two parts: A
geometric
part where certain local convexity data have to be verified, and a topological
part, a kind of ``Lokal-global-Prinzip'' \cite{HNP} which proves global
convexity \`a la Condevaux, Dazord, and Molino.

A further step was taken by Birtea, Ortega, and Ratiu \cite{BOR1, BOR2} who
consider a closed, not necessarily proper map $\mu\colon X\ra\widetilde{X}\ra
Y$, defined on a normal, first countable, arcwise connected Hausdorff space
$X$. The map $\mu$ has to be locally open onto its image, locally fiber
connected, having local convexity data. Using Va\v{\i}n\v{s}te\v{\i}n's lemma,
they prove that the light part
$\widetilde{X}\ra Y$ of $\mu$ is proper. This yields global convexity of
$\mu(X)$
for two almost disjoint kinds of target spaces $Y$, either the dual of a
Banach space \cite{BOR2} (which implies that the unit ball of $Y$ is weakly
compact), or a complete locally compact length metric space $Y$ \cite{BOR1}.
The second case applies to the cylinder-valued momentum map \cite{OR1, OR2},
another invention of Condevaux, Dazord, and Molino \cite{CDM}: For a
symplectic manifold $(X,\omega)$, the 2-form $\omega$ gives rise to a flat
connection on the trivial principal fiber bundle $X\times\mbox{Lie}(G)^\ast$
with holonomy group $H$.
The cylinder-valued momentum map $\overline{\mu}$ is obtained from $\mu$ by
factoring out $\overline{H}$ from the target space $Y$. The new target space
$\overline{\mu}(X)=Y/\overline{H}$ is a cylinder, hence geodesics on it may
differ from shortest paths. The convexity theorem then states that
$\overline{\mu}(X)$ is {\em weakly convex}, i.~e. any two points of
$\overline{\mu}(X)$ are connected by a geodesic arc.

In the present paper, we analyse the topological part of convexity, that is,
the passage from local to global convexity. We show that the
Lokal-global-Prinzip, as developed thus far, admits a substantial improvement
in at least three respects.

Firstly, we replace the monotone-light factorization $f\colon X\ra
\widetilde{X}\ra Y$ that was used for a momentum map $f=\mu$ by a new
factorization
$$f\colon X\stackrel{q^f}{\ltra} X^f\stackrel{f^{\#}}{\lra} Y$$
of any continuous map $f\colon X\ra Y$ which is locally open onto its image.
In a sense, $X^f$ is closer
to $Y$ than the leaf space $\widetilde{X}$ since $q^f\colon X\ra X^f$
factors through the monotone part $X\ra\widetilde{X}$ of $f$.
We show that $q^f$ is an open surjection, while $X^f$
admits a basis of open sets $U$ such that $f^{\#}$ maps $U$ homeomorphically
onto a subspace of $Y$ (Proposition~\ref{p5}).
Therefore, $f^{\#}$ can take the r\^ole of the light part of $f$, which means
that we can drop the assumption that $f$ (the momentum map) is locally fiber
connected.

Secondly, we concentrate on the target space $Y$ instead of $X$ to
derive the desired properties of $X^f$. In this way, the various assumptions
on $X$ boil down to a single one, namely, its connectedness as a topological
space. Nevertheless, we need no extra assumptions on the target space $Y$.

Thirdly, we merely assume that the map $f^{\#}$ is closed, a much weaker
condition than the closedness of $f$. Even the light part of $f$
need not be closed. For example, $f^{\#}$ is trivial for a local homeomorphism
$f$ - a light map which need not be closed, and with fibers of arbitrary
size. Using the properties of $Y$, we prove that the fibers of $f^{\#}$ are
finite (Proposition~\ref{p10}), so that the convexity structure of $Y$ can be
lifted along $f^{\#}$ (Theorem~\ref{t2}).

To make the interaction between convexity and topology more visible, we untie
the Lokal-global-Prinzip from its metric context by means of a general
concept of convexity, which might be of interest in itself. This also
unifies the two above mentioned types of target space considered in \cite{BOR1}
and \cite{BOR2}. In the linear case \cite{BOR2}, the target space $Y$ may be an
arbitrary (not necessarily complete) metrizable locally convex space instead
of a dual Banach space. (Metrizability is not needed unless the topology is
very strong, like in the case of a big locally convex direct sum.) In general,
geodesics in our target space $Y$ are one-dimensional continua which need not
be metrizable.

In previous versions of the Lokal-global-Prinzip, geodesic arcs or connecting
lines between two points of the target space $Y$ are obtained by
a metric on $Y$. Without a concept of length,
of course, geodesics are no longer available by shortening of arcs in the
spirit of the Hopf-Rinow theorem. Instead, we obtain geodesics by continued
{\em straightening}, using a local convexity structure. In other words, we
deal with a ``manifold'', that is, a Hausdorff space $Y$ covered by open
subspaces $U$ with an additional structure of convexity. The axioms of such
a {\em convexity space} $U$ are very simple: For any pair of points $x,y\in U$,
there is a minimal connected subset $C(x,y)$ containing $x$ and $y$, varying
continuously with the end points. In a topological vector space, $C(x,y)$ is
just the line segment between $x$ and $y$, while in a uniquely geodesic space,
$C(x,y)$ is the unique shortest path between $x$ and $y$. With respect to the
$C(x,y)$, there is a natural concept of convexity, and for a convexity space
$U$, we just require that the $C(x,y)$ are convex and that $U$ has a basis of
convex open sets (see Definition~\ref{d1}).

If convexity is given by a metric, straightening and shortening of
arcs leads to the same result, namely, a geodesic of minimal length.
For a non-metrizable arc $A$ between two points $x$ and $y$, there
is a substitute for the length of $A$, namely, the closed convex
hull $\overline{C(A)}$ which is diminished by straightening. As a
first step, an inscribed line path $L$ satisfies
$\overline{C(L)}\subs\overline{C(A)}$, and $\overline{C(L)}$ is the
closed convex hull of the finitely many extreme points of $L$. For a
given line path $L$ between $x$ and $y$, assume that the closed
convex hull $\overline{C(L)}$ is compact. Using Zorn's lemma, we
minimize the connected set $\overline{C(L)}$ to a compact convex set
$C$ with $x,y\in C$. In contrast to the Hopf-Rinow situation, where
the shortening of $L$ is achieved via the Arzel\`a-Ascoli theorem,
the straightening method needs the compactness of $\overline{C(L)}$
to show that connectedness carries over to $C$. By the local
convexity structure, it then follows that $C$ contains a line path
$L_0$ between $x$ and $y$. Thus if $C=L_0$, the line path $L_0$ must
be a geodesic.

So we require two properties to get the straightening process work:
First, the closed convex hull of a finite set must be compact; second, a
minimal compact connected convex set $C$ containing $x$ and $y$ has to be
a geodesic.

To establish a Lokal-global-Prinzip for continuous maps $X\ra Y$, possible
self-intersections of the arcs to be straightened have to be taken into
account. Precisely, this means that closed convex subsets of $Y$ have to be
replaced by {\em \'etale} maps, i.~e. closed locally convex maps
$e\colon C\ra Y$, such that the connected space $C$ admits a covering by open
sets mapped homeomorphically onto convex subsets of $Y$. We call
$Y$ a {\em geodesic manifold} if the above two properties hold with an adaption
to \'etale maps $e\colon C\ra Y$, that is, the second
property now states that if $C$ is compact and minimal with respect to
$x,y\in C$, then $e$ can be regarded as a geodesic with possible
self-intersections. (Such a geodesic is transversal if and only if $e=e^{\#}$.)
If the charts $U$ of $Y$ are
regular Hausdorff spaces which satisfy a finiteness condition (see
Definition~\ref{d2}) which holds, for example, if $U$ is either locally
compact or
first countable, we call $Y$ a {\em geodesic $q$-manifold}. Obvious examples
of geodesic $q$-manifolds are
complete locally compact length metric spaces, or metrizable locally convex
topological linear spaces (Examples~3 and 4).
Our main result consists in the following

\vspace{3mm}
\noindent {\bf Lokal-global-Prinzip.} \it Let $f\colon X\ra Y$ be a locally
convex continuous map from a connected topological space $X$ to a geodesic
$q$-manifold $Y$. Assume that $f^{\#}$ is closed. Then any two points of
$f(X)$ are connected by a geodesic arc.      \rm

For an inclusion map $f\colon C\hra Y$, the conditions on $f$ turn into the
assumptions of the Tietze-Nakajima theorem (see \cite{Nak}), i.~e. the subset
$C$ is closed, connected, and locally convex. Thus in case of a locally convex
topological vector space $Y$, the result for $C\hra Y$ yields Klee's convexity
theorem \cite{Klee}, while for a complete Riemannian manifold $Y$, we get a
theorem of Bangert \cite{Ban}.

\section{Convexity spaces}

Let $X$ be a Hausdorff space. We endow the power set $\PPP(X)$ with a topology
as follows. For any open set $U$ of $X$, define
\begin{equation}\label{1}
\widetilde{U}:=\{ C\in\PPP(X)\: |\: C\subs U\}.
\end{equation}
The collection $\BBB$ of sets (\ref{1}) is closed under finite intersection.
We take $\BBB$ as a basis of $\PPP(X)$.
\begin{Definition}\pt\label{d1}
\rm Let $X$ be a Hausdorff space together with a continuous map
\begin{equation}\label{2}
C\colon X\times X\ra\PPP(X).
\end{equation}
We call a subset $A\subs X$ {\em convex} if $C(x,y)\subs A$ holds for all
$x,y\in A$. We say that $X$ is a {\em convexity space} with respect to a map
(\ref{2}) if the following are satisfied.

\setlength{\partopsep}{-3mm}\begin{enumerate}
\setlength{\parskip}{-1mm}
\item[(C1)] The $C(x,y)$ are convex for all $x,y\in X$.
\item[(C2)] The $C(x,y)$ are minimal among the connected sets $C\subs X$ with
$x,y\in C$.
\item[(C3)] $X$ has a basis of convex open sets.
\end{enumerate}
\end{Definition}

Note that (C1) implies that $C(y,x)\subs C(x,y)$. Hence $C$ is symmetric:
\begin{equation}\label{3}
C(x,y)=C(y,x).
\end{equation}
From (C2) we infer that
\begin{equation}\label{4}
C(x,x)=\{ x\}.
\end{equation}
Moreover, (C2) implies that every convexity space $X$ is connected. The
restriction of the map (\ref{2}) to a convex subset $A\subs X$ makes $A$ into
a convexity space. Hence (C3) implies that $X$ is locally connected.
\begin{lem}\pt\label{l1}
Let $X$ be a convexity space. For $x,y\in X$, the set $C(x,y)\setm\{ y\}$ is
connected.
\end{lem}

\pf Let $A$ be the connected component of $x$ in $C(x,y)\setm\{ y\}$. Since
$\{ y\}$ is closed, every $z\in C(x,y)\setm\{ y\}$ admits a convex
neighbourhood $U$
with $y\notin U$. Hence $C(x,y)\setm\{ y\}$ is locally connected, and thus $A$
is open in $C(x,y)$. Since $C(x,y)$ is connected, it follows that $A$ cannot be
closed in $C(x,y)$. Thus $y\in\overline{A}$, which shows that $A\cup\{ y\}$
is connected. By (C2), this gives $A\cup\{ y\}=C(x,y)$, whence $A=C(x,y)\setm
\{ y\}$. \bx

As a consequence, the $C(x,y)$ can be equipped with a natural ordering.
\begin{prop}\pt\label{p1}
Let $X$ be a convexity space. For $x,y\in X$, the set $C(x,y)$ is linearly
ordered by
\begin{equation}\label{5}
z\le t\; :\Equ\; z\in C(x,t)\;\Equ\; t\in C(z,y)
\end{equation}
for $z,t\in C(x,y)$.
\end{prop}

\pf For any $z\in C(x,y)$, the set $C(x,z)\cup C(z,y)$ is connected. Therefore,
(C1) and (C2) give
\begin{equation}\label{6}
C(x,y)=C(x,z)\cup C(z,y).
\end{equation}
To verify the second equivalence in (\ref{5}), it suffices to show that
$$z\in C(x,t)\;\Ra\; t\in C(z,y)$$
holds for $z,t\in C(x,y)$. By Eq. (\ref{6}), it is enough to prove the
implication
\begin{equation}\label{7}
z\in C(x,t)\setm\{ t\}\;\Ra\; t\notin C(x,z).
\end{equation}
Assume that $z\in C(x,t)\setm\{ t\}$. Then Eq. (\ref{4}) gives $x\in C(x,t)
\setm\{ t\}$. Hence Lemma~\ref{l1}
and (C2) yield $C(x,z)\subs C(x,t)\setm\{ t\}$, which proves (\ref{7}).
Clearly, the relation (\ref{5}) is reflexive and transitive. By (\ref{7}), it
is a partial order. Furthermore, (\ref{5}) and (\ref{6}) imply that it is a
linear order. \bx

Note that the ordering of $C(x,y)$ depends on the pair $(x,y)$ which determines
the initial choice $x\le y$. Thus as an ordered set, $C(y,x)$ is dual to
$C(x,y)$.

\noindent {\bf Example 1.} Let $\Omega$ be a linearly ordered set. A subset $I$
of $\Omega$ is said to be an {\em interval} if $a\le c\le b$ with $a,b\in I$
implies that $c\in I$. The intervals $\{ c\in\Omega\:|\:c<b\}$ and
$\{ c\in\Omega\:|\:c>a\}$ with $a,b\in\Omega$ form a sub-basis for the {\em
order topology} of $\Omega$. Note that an open set of $\Omega$ is a
disjoint union
of open intervals. Therefore, $\Omega$ is connected if and only if it is a
{\em linear continuum}, i.~e. if every partition $\Omega=I\sqcup J$ into
non-empty intervals $I,J$ determines a unique element between $I$ and $J$.
With the order topology, a linear continuum $\Omega$ is a locally compact
convexity space with
\begin{equation}\label{8}
C(x,y)=\{ z\in\Omega\: |\: x\le z\le y\}
\end{equation}
in case that $x\le y$. Here the convex sets of $\Omega$ are just the connected
sets of $\Omega$.

\noindent {\bf Example 2.} More generally, we define a {\em tree continuum} to
be a Hausdorff space $X$ for which every two points $x,y\in X$ are contained in
a smallest connected set $C(x,y)$ such that each $C(x,y)$ is a linear
continuum, and $X$ carries the finest topology for which the inclusions
$C(x,y)\hra X$ are continuous. Thus $U\subs X$ is open if and only if every
$x\in U$ is an ``algebraically inner'' point (see \cite{Koe}, \S 16.2), i.~e.
if for each $y\in X\setm\{ x\}$, there exists some $z\in C(x,y)\setm\{x\}$
with $C(x,z)\setm\{z\}\subs U$. Then $X$ is a convexity space. For example,
every one-dimensional CW-complex without cycles is of this type.

\vspc
In the Euclidean plane $\R^2$, consider the solution curves $c\colon\R\ra\R^2$
of the differential equation $y'=3y^{\frac{3}{2}}$ (including the
singular solution $c\colon x\mapsto\tbinom{x}{0}$). With the finest topology
such that all solution curves are continuous, $\R^2$ becomes a tree continuum.
Here every point of the singular line is a branching point of order 4.

Note that a topological vector space $X$ is a convexity space with respect to
straight line segments if and only if $X$ is locally convex. The following
lemma is well-known (see \cite{Wil}, Theorem~26.15).
\begin{lem}\pt\label{l2}
Let $X$ be a connected topological space with an open covering $\UUU$. For any
pair of points $x,y\in X$, there is a finite sequence $U_1,\ldots,U_n\in\UUU$
with $x\in U_1$, $y\in U_n$, and $U_i\cap U_{i+1}\not=\varnothing$ for $i<n$.
\end{lem}

\begin{prop}\pt\label{p2}
Let $X$ be a convexity space. For $x,y\in X$, the subspace $C(x,y)$ is
compact and carries the order topology.
\end{prop}

\pf Let $C(x,y)=\bigcup\UUU$ be a covering by convex open sets. By
Lemma~\ref{l2}, there is a finite sequence $U_1,\ldots,U_n\in\UUU$
with $x\in U_1, y\in U_n$, and $U_i\cap U_{i+1}\not=\varnothing$ for $i<n$.
Hence $C(x,y)=U_1\cup\cdots\cup U_n$, which shows that $C(x,y)$ is compact.

\vspc
For $u<v$ in $C(x,y)$, the sets $C(x,u)$ and $C(v,y)$ are compact, hence
closed in $C(x,y)$. So the open intervals of $C(x,y)$ are open sets in
$C(x,y)$. Conversely, a convex open set in $C(x,y)$ is an interval
which must be an open interval since $C(x,y)$ is connected. \bx

Up to here, we have not used the continuity of the map (\ref{2}) in
Definition~\ref{d1}.
\begin{prop}\pt\label{p3}
Let $X$ be a convexity space. The closure of any convex set $A\subs X$ is
convex.
\end{prop}

\pf Let $A\subs X$ be a convex set, and let $x,y\in\overline{A}$ be given.
For any $z\in C(x,y)$, we
have to show that $z\in\overline{A}$. Suppose that there is a convex
neighbourhood $W$ of $z$ with $W\cap A=\varnothing$. Then $z\not=x,y$. By
Proposition~\ref{p2}, there exist $u,v\in W$ with $u<z<v$. Since $C(x,u)$
and $C(v,y)$ are compact, there are disjoint open sets $U,V$
in $X$ with $C(x,u)\subs U$ and $C(v,y)\subs V$ (see, e.~g., \cite{Kel},
chap.~V, Theorem~8). Hence $C(x,y)\subs U\cup V\cup W$. So
there are neighbourhoods $U'\subs U$ of $x$ and $V'\subs V$ of $y$ with
$C(x',y')\subs U\cup V\cup W$ for all $x'\in U'$ and $y'\in V'$. Choose $x',y'
\in A$. Then $C(x',y')\subs A$, which yields $C(x',y')\subs U\cup V$, where
$x'\in U'\subs U$ and $y'\in V'\subs V$, contrary to the connectedness of
$C(x',y')$. \bx

\begin{Definition}\pt\label{d2}
\rm Let $X$ be a convexity space. Define a {\em star} in $X$ with {\em center}
$x\in X$ and {\em end set} $E\subs X\setm\{ x\}$ to be a subspace $S(x,E):=
\bigcup\{ C(x,z)\:|\:z\in E\}$ with $C(x,z)\cap C(x,z')=\{x\}$ for
different $z,z'\in E$ such that
$S(x,E)$ carries the finest topology which makes the embeddings $C(x,z)\hra
S(x,E)$ continuous for all $z\in E$. We call $X$ {\em star-finite} if every
closed star in $X$ has a finite end set.
\end{Definition}
Thus every star is a tree continuum (Example~2). Recall that a topological
space $X$ is said to be a {\em $q$-space}
\cite{Mic1} if every point of $X$ has a sequence $(U_n)_{n\in\N}$ of
neighbourhoods such that every sequence $(x_n)_{n\in\N}$ with $x_n\in U_n$ has
an accumulation point. For example, every locally compact space, and every
first countable space $X$ is a $q$-space.

\begin{prop}\pt\label{p4}
Let $X$ be a convexity space which is a $q$-space. Then $X$ is star-finite.
\end{prop}

\pf Let $S(x,E)$ be a closed star in $X$, and let $(U_n)_{n\in\N}$ be
a sequence of neighbourhoods of $x$ such that every sequence
$(x_n)_{n\in\N}$ with $x_n\in U_n$ has an accumulation point. Suppose that $E$
is infinite. Since $U_n\cap C(x,z)\not=\{x\}$ for all $n\in\N$ and $z\in E$,
we find a subset $\{z_n\: |\:n\in\N\}$ of $E$ and a sequence $(x_n)_{n\in\N}$
with $x\not=x_n\in C(x,z_n)\cap U_n$. Thus $(x_n)_{n\in\N}$ has an
accumulation point $z$. Because of the star-topology, $z$ cannot belong to
$S(x,E)$, contrary to the assumption that $S(x,E)$ is closed. \bx

\section{Local openness onto the image}

For a topological space $X$, the infinitesimal structure at a point $x$ is
given by the set $\DDD_x$ of filters on $X$ which converge to $x$. Let
$\FFF(X)$ denote the set of all filters on $X$. We make $\FFF(X)$ into a
topological space with a basis of open sets
\begin{equation}\label{9}
\widetilde{U}:=\{\alpha\in\FFF(X)\: |\: U\in\alpha\},
\end{equation}
where $U$ runs through the class of open sets in $X$. Every continuous map
$f\colon X\ra Y$ induces a map $\FFF(f)\colon\FFF(X)\ra\FFF(Y)$. For an open
set $V$ in $Y$, we have
\begin{equation}\label{10}
\FFF(f)^{-1}(\widetilde{V})=\widetilde{f^{-1}(V)},
\end{equation}
which shows that $\FFF(f)$ is continuous.
Consider the subspace
\begin{equation}\label{11}
\DDD(X):=\{(x,\alpha)\in X\times\FFF(X)\: |\:\alpha\in\DDD_x\}
\end{equation}
of $X\times\FFF(X)$.
Note that for every $x\in X$, the neighbourhood filter $\UU(x)$ of $x$ is the
coarsest filter in $\DDD_x$. Thus, regarding $\DDD_x$ as a subset of $\DDD(X)$,
we get a pair of continuous maps
\begin{equation}\label{12}
\begin{diagram}[midshaft]
X & \rTo^{\mbox{\scriptsize $\UU$}} & \DDD(X) & \rOn^{\mbox{\scriptsize lim}}
& X
\end{diagram}
\end{equation}
with $\mbox{lim}(x,\alpha):=x$ and $\mbox{lim}\circ\UU=1_X$. In particular,
$\DDD_x=\mbox{lim}^{-1}(x)$.

For a continuous map $f\colon X\ra Y$, the local behaviour at
$x\in X$ is given by the induced map $\DDD_x f\colon\DDD_x\ra
\DDD_{f(x)}$. Thus we get an endofunctor $\DDD\colon\mbox{\bf Top}\ra
\mbox{\bf Top}$ of the category \mbox{\bf Top} of topological spaces with
continuous maps as morphisms. The functor $\DDD$ is augmented by the natural
transformation $\mbox{lim}\colon\DDD\ra 1$.

\begin{Definition}\pt\label{d3}
\rm A continuous map $f\colon X\ra Y$ between topological spaces is said to
be {\em locally open onto its image} \cite{Ben} if every $x\in X$ admits an
open neighbourhood $U$ such that the induced map $U\tra f(U)$ is open onto the
subspace $f(U)$ of $Y$. We call $f$ {\em filtered} if $f$ is locally open onto
its image and $\DDD(f)\circ\UU$ is injective.
\end{Definition}

We have the following structure theorem for continuous maps which are locally
open onto its image.
\begin{prop}\pt\label{p5}
Let $f\colon X\ra Y$ be a continuous map which is locally open onto its image.
Up to isomorphism, there is a unique factorization $f=pq$ into an open
surjection $q$ and a
filtered map $p$. If $f$ is filtered, then every point $x\in X$ has an open
neighbourhood which is mapped homeomorphically onto a subspace of $Y$.
\end{prop}

\pf Consider the following commutative diagram
$$\begin{diagram}[s=6ex,midshaft]
1\colon\hspace*{-10mm} & X & \rTo^{\mbox{$\UU$}} & \DDD(X) & \rTo^{\mbox{lim}}
& X\\
& \dOn>{\mbox{$q^f$}} & & \dTo>{\mbox{$\DDD(f)$}} & & \dTo>{\mbox{$f$}}\\
f^{\#}\colon\hspace*{-6mm} & X^f & \rInc^{\mbox{$e$}} & \DDD(Y) &
\rTo^{\mbox{lim}} & Y\dpc\\
\end{diagram} $$
where $X^f$ is the image of $\DDD(f)\circ\UU$, regarded as a quotient space of
$X$, and $f^{\#}:=\mbox{lim}\circ e$.
We will prove that $f=f^{\#}\circ q^f$ gives the desired factorization. Let us
show first that $q^f$ is open. Thus let $U$ be an open set of $X$. We have to
verify that $(q^f)^{-1}q^f(U)$ is open in $X$. Since $f$ is locally open onto
its image, we can assume that the induced map $U\tra f(U)$ is
open. Let $x\in (q^f)^{-1}q^f(U)$ be given. Then $q^f(x)\in q^f(U)$. So there
exists some $y\in U$ with $q^f(x)=q^f(y)$, i.~e. $f(x)=f(y)$ and
$f(\UU(x))=f(\UU(y))$. Hence there is an open neighbourhood $V\in\UU(x)$ with
$f(V)\subs f(U)$. Again, we can assume that the induced map $V\tra f(V)$ is
open. Furthermore, there is an open neighbourhood $U'\subs U$ of $y$ with
$f(U')\subs f(V)$, and $f(U')$ is open in $f(U)$, hence in $f(V)$. Therefore,
$V':=V\cap f^{-1}(f(U'))$ is an open neighbourhood of $x$ with $f(V')=f(U')$.

\vspc
For any $x'\in V'$, there is a point $y'\in U'$ with $f(x')=f(y')$. So the
continuity of $f$ implies that $f(\UU(x'))=f(\UU(y'))$, which gives
$q^f(x')=q^f(y')$, and thus $V'\subs (q^f)^{-1}q^f(U')\subs (q^f)^{-1}q^f(U)$.
This proves that
$q^f$ is open. Consequently, $f^{\#}$ is locally open onto its image.

\vspc
Since $q^f$ is open, we have a commutative diagram
$$\begin{diagram}[w=8ex,h=5ex,midshaft]
X & \rOn^{\mbox{$q^f$}} & X^f\\
\dTo>{\mbox{$\UU$}} & & \dTo>{\mbox{$\UU$}}\\
\DDD(X) & \rTo^{\mbox{$\DDD(q^f)$}} & \DDD(X^f).\\
\end{diagram} $$
Hence $\DDD(f^{\#})\circ\UU\circ q^f=\DDD(f^{\#})\circ\DDD(q^f)
\circ\UU=\DDD(f)\circ\UU=e\circ q^f$. Therefore, $\DDD(f^{\#})\circ\UU=e$,
which implies that $f^{\#}$ is filtered.

\vspc
Now let $f=pq=p'q'$ be two factorizations with $p,p'$ filtered and $q,q'$
open. Then $\DDD(p')\circ\UU\circ q'=\DDD(p')\circ\DDD(q')\circ\UU=\DDD(p)
\circ\DDD(q)\circ\UU=\DDD(p)\circ\UU\circ q$. Since $\DDD(p')\circ\UU$ is
injective, there exists a continuous map $e\colon E\ra E'$ such that $q'=eq$.
So we get a commutative diagram
$$\begin{diagram}[s=5ex]
X & \rOn^{\mbox{$q$}} & E & \rTo^{\mbox{$p$}} & Y\\
\deq & & \dTo>{\mbox{$e$}} & & \deq\\
X & \rOn^{\mbox{$q'$}} & E' & \rTo^{\mbox{$p'$}} & Y\\
\end{diagram} $$
By symmetry, we find a continuous map $e'\colon E'\ra E$ with $q=e'q'$ and
$p'=pe'$. Since $q$ and $q'$ are surjective, $e$ must be a homeomorphism. This
proves the uniqueness of the factorization.

\vspc
Finally, let $f\colon X\ra Y$ be filtered. For a given $x\in X$, let $U$ be an
open neighbourhood such that the induced map $r\colon U\tra f(U)$ is open.
Since $i\colon U\hra X$ is open, we have a commutative diagram
$$\begin{diagram}[w=8ex]
& & X & \rTo^{\mbox{$\UU$}} & \DDD(X) & & \\
& \ldTo \makebox[0pt][l]{\hspace*{-2.5mm}\raisebox{2mm}{$f$}} &
\uInc>{\mbox{$i$}} & & \uTo>{\mbox{$\DDD(i)$}} & \rdTo
\makebox[0pt][l]{\hspace*{-.5mm}\raisebox{1.5mm}{$\DDD(f)$}} & \\
\mbox{\hspace*{4.9mm} $Y$\rule[-.77mm]{0pt}{4.1mm}} & & U &
\rTo^{\mbox{$\UU$}} & \DDD(U) & & \DDD(Y) \\
& \luInc \makebox[0pt][l]{\hspace*{-2.5mm}\raisebox{-2mm}{$j$}} &
\dOn>{\mbox{$r$}} & & \dTo>{\mbox{$\DDD(r)$}} & \ruTo
\makebox[0pt][l]{\hspace*{-.5mm}\raisebox{-2mm}{$\DDD(j)$}} & \\
& & f(U) & \rTo^{\mbox{$\UU$}} & \DDD(f(U)) & & \\
\end{diagram} $$
which shows that $\DDD(j)\circ\UU\circ r=\DDD(f)\circ\UU\circ i$ is injective.
Hence $r$ is injective. \bx

In the sequel, we keep the notation of Proposition~\ref{p5} and write
\begin{equation}\label{13}
f\colon X\stackrel{q^f}{\ltra} X^f\stackrel{f^{\#}}{\lra} Y
\end{equation}
for the factorization of a map $f$ which is locally open onto its image.

\noindent {\bf Remarks. 1.} Although the factorization (\ref{13}) is unique
up to
isomorphism, it does not give rise to a factorization system \cite{FK, CHK},
i.~e. a pair $(\EE,\MM)$ of subcategories such that every commutative square
\begin{equation}\label{14}
\begin{diagram}
E_1 & \rTo^{\mbox{$f_1$}} & M_1\\
\dTo<{\mbox{$e$}} & \makebox[0pt][l]{\hspace*{-5.5mm} $d$}\rudash &
\dTo>{\mbox{$m$}}\\
E_0 & \rTo^{\mbox{$f_0$}} & M_0\\
\end{diagram}
\end{equation}
with $e\in\EE$ and $m\in\MM$ admits a unique diagonal $d$ with $f_1=de$ and
$f_0=md$ (see \cite{Her},
Proposition~1.4). Apart from the fact that local openness onto the image is not
closed under composition (consider the maps $\R\stackrel{i}{\hra}\R^2
\stackrel{p}{\tra}\R$ with $i(x)=\binom{x}{x^3-3x}$ and $p\colon\binom{x}{y}
\mapsto y$), there cannot be a factorization system since open surjections are
not stable under pushout (take, e.~g., the pushout of the open surjection
$\R\tra\{0\}$ and the inclusion $\R\hra\R^2$).

{\bf 2.} If $f\colon X\ra Y$ is locally open onto its image and locally fiber
connected \cite{Ben, HNP}, the lemma of Benoist (\cite{Ben}, Lemma~3.7)
states that the monotone part $\pi$ of the monotone-light factorization
$f=\widetilde{f}\circ\pi$ is open. Here the local fiber-connectedness of $f$
implies that $\pi$ is locally open onto its image. Hence $\pi=q^\pi$
is open by Proposition~\ref{p5}. In general, $q^f$ always factors through
$\pi$, but the two factorizations need not be isomorphic. For example, a local
homeomorphism $f\colon X\tra Y$ is open, but its fibers are discrete.

\section{Convexity of maps}

In this brief section, we introduce local convexity and extend this concept
from subsets to continuous maps (cf. \cite{KB} for a notion of convex maps
in terms of paths).
\begin{Definition}\pt\label{d4}
\rm Let $X$ be a topological space. We define a {\em local convexity
structure}
on $X$ to be an open covering $X=\bigcup\UUU$ by convexity spaces $U\in\UUU$
(with the induced topology) such that for any $U\in\UUU$, every convex open
subspace of $U$ belongs to $\UUU$ (as a convexity space). We call a subset
$C\subs X$ {\em convex} if $C\cap U$ is convex for all $U\in\UUU$. We say
that $C$ is {\em locally convex} if every $z\in C$ admits a neighbourhood
$U\in\UUU$ such that $C\cap U$ is convex.
\end{Definition}
The covering $\UUU$ will be referred to as the {\em atlas} of the
local convexity structure. In the special case $X\in\UUU$, the atlas $\UUU$
just consists of the convex open sets of a convexity space $X$.

In contrast to local convexity, our concept of convexity refers to all
sets in $\UUU$. So the
intersection of convex sets is convex, and every subset $A\subs X$ admits a
{\em convex hull} $C(A)$, that is, a smallest convex set $C\sups A$. The next
proposition generalizes Proposition~\ref{p3}.
\begin{prop}\pt\label{p6}
Let $X$ be a topological space with a local convexity structure $\UUU$. The
closure of any convex set $A\subs X$ is convex.
\end{prop}

\pf For every $U\in\UUU$, we have $\overline{A}\cap U=\overline{A\cap U}
\cap U$. This set is convex by Proposition~\ref{p3}. Hence $A$ is convex. \bx

Definition~\ref{d4} admits a natural extension to continuous maps.

\begin{Definition}\pt\label{d5}
\rm Let $f\colon X\ra Y$ be a continuous map between topological spaces, where
$Y$ has a local convexity structure $\VVV$. We call $f$
{\em locally convex} if every $x\in X$ admits an open neighbourhood $U$ such
that the induced map $U\tra f(U)$ is open, and $f(U)$ is a convex subspace
of some $V\in\VVV$.
\end{Definition}

\noindent {\bf Remarks. 1.} A subset $A\subs Y$ is locally convex if and only
if the inclusion map $A\hra Y$ is locally convex.

\vspc
{\bf 2.} The open neighbourhood $U$ of $x$ in Definition~\ref{d5} can be chosen
arbitrarily small. In fact, let $U'\subs U$ be any smaller open neighbourhood
of $x$. Then $f(U')$ is an open subset of $f(U)$. Hence there exists some
$V'\in\VVV$ with $f(x)\in V'\cap f(U)\subs f(U')$. Thus $U'':=
U'\cap f^{-1}(V')$ is an open neighbourhood of $x$ with $f(U'')=V'\cap f(U')=
V'\cap f(U)$, which is a convex subspace of $V'$.

\vspc
{\bf 3.} If $X$ is a connected Hausdorff space and $Y$ a length metric space
\cite{BrH, Gr}, a continuous map $f\colon X\ra Y$ is locally convex if and only
if $f$ is locally open onto its image and has local convexity data in the
sense of \cite{BOR1}.

\begin{prop}\pt\label{p7}
Let $f\colon X\ra Y$ be a continuous map between topological spaces, where $Y$
has a local convexity structure $\VVV$. If $f$ is
locally convex, then $f^{\#}$ is locally convex.
\end{prop}

\pf Assume that $f$ is locally convex, and let $U$ be an open neighbourhood of
$x\in X$ such that the induced map $U\tra f(U)$ is open onto a convex
subspace of some $V\in\VVV$. Since $q^f$ is open by
Proposition~\ref{p5}, this property of $U$ carries over to the neighbourhood
$q^f(U)$ of $q^f(x)$. Hence $f^{\#}$ is locally convex. \bx

\section{Geodesic manifolds}

In this section, we introduce a general concept of geodesic which does not
refer to any kind of metric.

\begin{Definition}\pt\label{d6}
\rm Let $Y$ be a topological space with a local convexity structure $\VVV$,
and let $e\colon C\ra Y$ be a continuous map with a connected topological space
$C$. By $\VVV_e$ we denote the set of all open sets $U$ in $C$ which are mapped
homeomorphically onto a convex subspace of some $V\in\VVV$. We call $e$
{\em \'etale} if $e$ is closed and $\VVV_e$ covers $C$. We say that
$e\colon C\ra Y$ is {\em generated} by a subset $F\subs C$ if
there is no closed connected subspace $A\subsn C$ with $F\subs A$ such that
$e(U\cap A)$ is convex for all $U\in\VVV_e$.
\end{Definition}
In particular, \'etale maps are locally convex. Furthermore,
every \'etale map $e\colon C\ra Y$ induces a
local convexity structure $\VVV_e$ on $C$. If $F\subs C$ is connected, then
$C(F)$ is connected, and an
\'etale map $e\colon C\ra Y$ is generated by $F$ if and only if
$\overline{C(F)}=C$. Note that the composition of \'etale maps is \'etale.

\begin{Definition}\pt\label{d7}
\rm Let $Y$ be a Hausdorff space with a local convexity structure $\VVV$. We
call $Y$ a {\em geodesic manifold} if the following are satisfied.

\setlength{\partopsep}{-3mm}\begin{enumerate}
\setlength{\parskip}{-1mm}
\item[(G1)] For a finite set $F\subs Y$, the closure of $C(F)$ is compact.
\item[(G2)] If an \'etale map $e\colon C\ra Y$ with $C$ compact is generated by
$\{x,y\}\subs C$, then every connected set $A\subs C$ with $x,y\in A$
coincides with $C$.
\end{enumerate}
If, in addition, the $V\in\VVV$ are star-finite and regular (as topological
spaces), we call $Y$ a {\em geodesic $q$-manifold}.
\end{Definition}
The letter ``$q$'' is reminiscent of Proposition~\ref{p4}. Since a geodesic
manifold $Y$ is locally connected, \cite{Bou}, chap.~I, 11.6,
Proposition~11, implies that $Y$ is the topological sum of its connected
components.

\begin{Definition}\pt\label{d8}
\rm Let $Y$ be a geodesic manifold. We define a {\em geodesic} in $Y$ to be
an \'etale map $e\colon C\ra Y$, generated by $\{ x,y\}\subs C$, where $C$ is
compact. The points
$e(x)$ and $e(y)$ will be called the {\em end points} of the geodesic.
\end{Definition}

More generally, we define a {\em line path} in $Y$ to be a continuous map
\mbox{$e\colon L\ra Y$}, where $L$ is a linear continuum (Example~1) with end
points $x_0$ and $x_n$ and a sequence of intermediate points $x_0<x_1<\cdots<
x_n$ such that for $i<n$,
the restriction of $e$ to the interval $[x_i,x_{i+1}]$ is an inclusion which
identifies $[x_i,x_{i+1}]$ with $C(e(x_i),e(x_{i+1}))\subs U_i$ for some $U_i$
in the atlas of $Y$. If $e$ is an inclusion, we speak of a {\em simple} line
path and identify it with the subset $L\subs Y$. A subset $A\subs Y$ will be
called {\em line-connected} if every pair of points
$x,y\in A$ is connected by a simple line path $L\subs A$.

\begin{prop}\pt\label{p8}
Let $Y$ be a geodesic manifold with atlas $\VVV$, and let $e\colon C\ra Y$
be an \'etale map. Then $C$ is line-connected.
\end{prop}

\pf Let $x,y\in C$ be given. By Lemma~\ref{l2}, there is a sequence $U_1,
\ldots,U_n\in\VVV_e$ with $x\in U_1, y\in U_n$, and $U_i\cap U_{i+1}\not=
\varnothing$ for $i<n$. Choose $x_i\in U_i\cap U_{i+1}$ for $i<n$. With $x_0:=
x$ and $x_n:=y$, the $C(x_i,x_{i+1})$ constitute a line path $e\colon L\ra Y$
in
$C$ which connects $x$ and $y$. Assume that the interval $[x,x_i]\subs L$ maps
onto a simple line path $L'$. If $C(x_i,x_{i+1})$ intersects $L'$ in a point
$\not=x_i$, there is a largest $z\in C(x_i,x_{i+1})$ with property. Thus, if
$z'$ denotes the corresponding point on $L'$, we can replace the interval
$[z',z]$ by $\{z\}$ and attach the segment $C(z,x_{i+1})$. After finitely many
steps, we obtain a simple line path between $x$ and $y$. \bx

By (G2), we have the following

\vspace{3mm}
\noindent {\bf Corollary.} \it Let $Y$ be a geodesic manifold. Every geodesic
with end points $x,y\in Y$ is a line path.  \rm

\vspace{3mm}
In particular, a simple geodesic with end points $x,y\in Y$ is just a minimal
connected set $C\subs Y$ with $x,y\in C$ which is locally convex.

Let $Y$ be a geodesic manifold. For $x,y\in Y$, we define a {\em simple arc}
between $x$ and $y$ to be a subspace $A\subs Y$ which is a linear continuum
with end points $x$ and $y$. We fix a
linear order on $A$ such that $x$ becomes the smallest element and denote the
set of all such $A$ by $\Omega_Y(x,y)$. In particular, every simple line path
between $x$ and $y$ belongs to $\Omega_Y(x,y)$. Clearly, every $A\in
\Omega_Y(x,y)$ admits an inscribed line path $L$ between $x$ and $y$.
Although there is no concept of length at our disposal, the intuition that
$L$ is ``shorter'' than $A$ can be expressed by the inclusion $\overline{C(L)}
\subs\overline{C(A)}$. Thus it is natural to define a preordering on
$\Omega_Y(x,y)$ by
\begin{equation}\label{15}
A\prec B :\Equ\;\overline{C(A)}\subs\overline{C(B)}.
\end{equation}
If $A\prec B$ holds for a pair $A,B\in
\Omega_Y(x,y)$, we say
that $A$ is a {\em straightening} of $B$.
Define $B\in\Omega_Y(x,y)$ to be {\em minimal} if $A\prec B$ implies
$B\prec A$ for all $A\in\Omega_Y(x,y)$. We have the following straightening
theorem which justifies the term ``geodesic'' manifold in Definition~\ref{d7}.
\begin{thm}\pt\label{t1}
Let $Y$ be a geodesic manifold. Every simple arc $A\in\Omega_Y(x,y)$ in $Y$
can be straightened to a minimal $C\in\Omega_Y(x,y)$. A simple arc
$A\in\Omega_Y(x,y)$ is minimal if and only if $A$ is a convex simple geodesic.
\end{thm}

\pf Let $A\in\Omega_Y(x,y)$ be given. Since $C(A)$ is connected,
$\overline{C(A)}$ is connected. Proposition~\ref{p6} implies that
$\overline{C(A)}$ is convex. So the inclusion $\overline{C(A)}\hra Y$ is
\'etale. By
Proposition~\ref{p8}, there exists a simple line path $L\subs\overline{C(A)}$
between $x$ and $y$. Hence $L\prec A$. As $L$ belongs to the convex hull of a
finite set, (G1) implies that $\overline{C(L)}$ is compact. We have to
verify that $\overline{C(L)}$ contains a minimal
$C\in\Omega_Y(x,y)$. Let $\CC$ be a chain of compact convex connected sets
$C\subs\overline{C(L)}$ with $x,y\in C$. Then $D:=\bigcap\CC$ is compact and
convex, and $x,y\in D$. We show first that every open set $V$ of $Y$ with
$D\subs V$ contains some $C\in\CC$. In fact, the set $\overline{C(L)}$
is compact, and $\bigcap_{C\in\CC}(C\setm V)=\varnothing$. Hence $C\setm V=
\varnothing$ for some $C\in\CC$.

\vspc
Next we show that $D$ is connected. Suppose that there is a disjoint union
$D=D_1\sqcup D_2$ with non-empty compact sets $D_1$ and $D_2$. Then we can
find open sets $U_1$ and $U_2$ in $Y$ with $D_i\subs U_i$ such that $U_1\cap
U_2=\varnothing$ (see, e.~g., \cite{Kel}, chap. V, Theorem~8).
Hence $D\subs U_1\sqcup U_2$, which yields $C\subs U_1\sqcup U_2$ for some
$C\in\CC$. Since $C$ is connected, we can assume that $C\subs U_1$.
This gives $D_2\subs U_1\cap U_2=\varnothing$, a contradiction. Thus $D$ is a
connected. By Zorn's lemma, it follows that there exists a minimal compact
convex connected set $C$ with $x,y\in C$. Hence $C\hra Y$ is an \'etale
map generated by $\{ x,y\}$. Therefore, (G2) implies that $C$ admits no
connected proper subset $C'\subs C$ with $x,y\in C'$. By Proposition~\ref{p8},
it follows that $C$ is a simple line path, whence $C\in\Omega_Y(x,y)$, and
$C$ is minimal.

\vspc
In particular, we have shown that if $A\in\Omega_Y(x,y)$ is minimal, then
$A$ is a convex simple geodesic between $x$ and
$y$. Conversely, if $A\in\Omega_Y(x,y)$ is a convex simple geodesic, then
$A=\overline{C(A)}$, and thus $A$ is minimal. \bx

We conclude this section with some typical examples.

\noindent {\bf Example 3.} Let $Y$ be a geodesic manifold with atlas $\VVV$,
and let $Z$ be a closed locally convex subspace. Then $Z\hra Y$ is \'etale.
Every finite set $F$ in $Z$ is contained in a compact convex set $C$ in $Y$.
Hence $C\cap Z$ is compact and convex in $Z$. Thus $Z$ satisfies (G1). As (G2)
trivially carries over to $Z$, it follows that $Z$ is a geodesic manifold.
If $Y$ is a geodesic $q$-manifold, then so is $Z$.

\noindent {\bf Example 4.} Let $Y$ be a complete locally compact length metric
space \cite{BrH, Gr}. By the Hopf-Rinow theorem (\cite{BrH},
Proposition~I.3.7),
the closed metric balls in $Y$ are compact, and any two points in $Y$ are
connected by a shortest path. It is natural to assume that $Y$ admits a basis
of convex open sets where shortest paths are unique. This provides $Y$ with a
local convexity structure $\VVV$ which satisfies (G1). Note that by
\cite{BrH}, I.3.12, the map (\ref{2}) is continuous where it is defined.

\vspc
Now let $e\colon C\ra Y$ be an \'etale map generated by $\{x,y\}\subs C$,
where $C$ is compact. Similar to the case of a covering of length metric
spaces (\cite{BrH}, Proposition~I.3.25), the length metric $d_Y$ of $Y$ can be
lifted to a length metric $d_C$ of $C$ such that $d_C(u,v)\ge d_Y(e(u),e(v))$
for all $u,v\in C$. (If $d_C(u,v)=0$ with $u\not=v$, a neighbourhood
$U\in\VVV_e$ of $u$ cannot contain $v$. As $U$ contains a closed
neighbourhood of $u$ in $C$, we get $d_C(u,v)>0$.) Since $C$ is compact, the
Hopf-Rinow theorem, applied to $C$, yields a shortest path $L\subs C$ between
$x$ and $y$. Hence $C=L$, which proves (G2). By Proposition~\ref{p4}, $Y$ is a
geodesic $q$-manifold.

\noindent {\bf Example 5.} Let $Y$ be a locally convex topological vector
space. For $x,y\in Y$, we set $C(x,y):=\{ \lambda x+(1-\lambda)y\: |\:
\lambda\ge 0\}$ to make $Y$ into a convexity space. For a finite set
$F\subs Y$, the closed convex hull $\overline{C(F)}$ of $F$ is contained in a
finite dimensional subspace of $Y$. Hence $\overline{C(F)}$ is compact. Thus
$Y$ satisfies (G1). Let
$e\colon C\ra Y$ be an \'etale map generated by $\{x,y\}\subs C$,
where $C$ is compact. By Proposition~\ref{p8}, $e$ is generated by a
simple line path in $C$. Hence $e(C)$ is contained in a finite dimensional
subspace of $Y$. So Example~4 applies, which proves (G2). Thus $Y$ is a
geodesic manifold. If $Y$ is metrizable, i.~e. first countable (\cite{Sch},
I, Theorem~6.1), then $Y$ is a geodesic $q$-manifold by Proposition~\ref{p4}.

\section{The Lokal-global-Prinzip}

With respect to convex neighbourhoods, \'etale maps have the following
disjointness property.
\begin{prop}\pt\label{p9}
Let $Y$ be a geodesic manifold with atlas $\VVV$, and let $e\colon C\ra Y$
be an \'etale map. Assume that $U,U'\in\VVV_e$.
If $e|_{U\cup U'}$ is not injective, then $U\cap U'=\varnothing$.
\end{prop}

\pf If $e|_{U\cup V}$ is not injective, there exist $x\in U$ and $x'\in U'$
with $e(x)=e(x')$. Suppose that there is some $z\in U\cap U'$. Then $x\not=z$,
and $U\cap U'\cap C(x,z)$ is a convex open subset of $C(x,z)\setm\{ x\}$.
Hence there
is a point $t\in C(x,z)$ with $(U\setm U')\cap C(x,z)=C(x,t)$. So the
homeomorphisms $C(x,z)\cong C(e(x),e(z))\cong C(x',z)$ give rise to a point
$t'\in U'$ with $e(t)=e(t')$ and $(U'\setm U)\cap C(x',z)=C(x',t')$. Moreover,
$D:=C(t,z)\cup C(t',z)=C(t,z)\cup\{ t'\}$. Therefore, $D$ is not a minimally
connected superset of $\{t,z\}$. On the other hand, $D$ is compact with open
subsets $C(t,z)$ and $C(t',z)$. Hence $e|_D\colon D\ra Y$ is an \'etale map
generated by $\{t,z\}$, contrary to (G2). \bx

As an immediate consequence, the fibers of an \'etale map can be separated by
pairwise disjoint neighbourhoods.

\vspace{3mm}
\noindent {\bf Corollary 1.} \it Let $Y$ be a geodesic manifold, and let
$e\colon C\ra Y$ be an \'etale map. For a given $y\in Y$,
choose a neighbourhood $U_x\in\VVV_e$ of each $x\in f^{-1}(y)$. Then the $U_x$
are pairwise disjoint.        \rm

\vspace{3mm}
\noindent {\bf Corollary 2.} \it Let $Y$ be a geodesic manifold, and let
$e\colon C\ra Y$ be an \'etale map. Then $C$ is a Hausdorff space. \rm

\vspace{3mm}
\pf Let $x,x'\in C$
be given. If $e(x)\not= e(x')$, there are disjoint neighbourhoods of $e(x)$
and $e(x')$, and their inverse images give disjoint neighbourhoods of $x$ and
$x'$. So we can assume that $e(x)=e(x')$. Choose $U,U'\in\VVV_e$ with $x\in U$
and $x'\in U'$. By Proposition~\ref{p9}, $U\cap U'=\varnothing$. Thus $C$ is
Hausdorff. \bx

If the geodesic manifold is regular, the fibers are even discrete, which leads
to the following finiteness result.

\begin{prop}\pt\label{p10}
Let $e\colon C\ra Y$ be an \'etale map into a geodesic $q$-manifold $Y$. Then
the fibers of $e$ are finite.
\end{prop}

\pf Let $\VVV$ denote the atlas of $Y$, and let $y\in Y$ be given. For each
$x\in e^{-1}(y)$, we choose a neighbourhood $U_x\in\VVV_e$ such that the
images $e(U_x)$ are contained in a fixed $V'\in\VVV$. By the
Corollary~1, these neighbourhoods are pairwise disjoint. Without loss of
generality, we can assume that $|C|>1$. Since $C$ is a connected Hausdorff
space by Corollary~2, this implies that $C$ has no isolated points. As $e$ is
closed, the complement of $\bigcup\{ U_x\: |\: x\in e^{-1}(y)\}$ is mapped to
a closed set $A\subs Y$ with $y\notin A$. So there exists an open
neighbourhood $W\subs V'$ of $y$ with $e^{-1}(W)\subs\bigcup\{ U_x\: |\: x\in
e^{-1}(y)\}$. By the regularity of $Y$, we find a convex open neighbourhood
$V$ of $y$ with $\overline{V}\subs W$.

\vspc
For any $x\in e^{-1}(y)$, the set $U_x\cap e^{-1}(V)$ is an open neighbourhood
of $x$, hence not a singleton. Therefore, the $V_x:=e(U_x\cap e^{-1}(V))$
are convex subsets of $V$ with $|V_x|>1$ and $y\in V_x$. Choose arbitrary
$z_x\in U_x\cap e^{-1}(V)$ with $y_x:=e(z_x)\not=y$ for all $x\in e^{-1}(y)$.
Now let $Z\subs\bigcup\{C(x,z_x)\:|\:x\in e^{-1}(y)\}$ be such that
$Z\cap C(x,z_x)$ is closed in $U_x\cap e^{-1}(V)$ for every $x\in e^{-1}(y)$.
We claim that $Z$ is closed. Thus let $z\in\overline{Z}$ be given. Then
$e(z)\in\overline{e(Z)}\subs\overline{V}\subs W$. Hence
$z\in e^{-1}(W)\subs\bigcup\{U_x\: |\:\mbox{$x\in e^{-1}(y)$}\}$, which
yields $z\in Z$. Thus $Z$ is closed. Since $e$ is closed, this implies
that $S(y):=\bigcup\{ C(y,y_x)\: |\: x\in e^{-1}(y)\}$ is closed and carries
the finest topology such that the maps $C(y,y_x)\hra S(y)$ are continuous for
all $x\in e^{-1}(y)$.

\vspc
Suppose that $e^{-1}(y)$ is infinite. Since $S(y)$ cannot be a closed star,
there must be an infinite countable subset $E$ of $e^{-1}(y)$ with
$C(y,y_u)\cap C(y,y_v)\not=\{y\}$ for different $u,v\in E$. Hence there is a
point $y'\in V\setm\{y\}$ and a set $Z\subs\bigcup\{C(x,z_x)\:|\:x\in
e^{-1}(y)\}$ with $|Z\cap C(x,z_x)|=1$ for all $x\in E$
such that $e(Z)$ is an infinite non-closed subset of $C(y,y')$. Since $Z$ is
closed, this gives a contradiction. \bx

As a consequence, the geodesic structure of a geodesic $q$-manifold can be
lifted along \'etale maps.

\begin{thm}\pt\label{t2}
Let $e\colon C\ra Y$ be an \'etale map into a geodesic $q$-manifold $Y$ with
atlas $\VVV$. Then $C$ is a geodesic $q$-manifold with atlas $\VVV_e$.
\end{thm}

\pf By Corollary~2 of Proposition~\ref{p9}, $C$ is a Hausdorff space. We show
first that $C$ is regular. Let $U_x\in\VVV_e$ be a neighbourhood
of $x\in C$. We choose neighbourhoods $U_z\in\VVV_e$ for all $z$ in the fiber
of $y:=e(x)$. By Corollary~1 of Proposition~\ref{p9}, the $U_z$ are pairwise
disjoint. Since $Y$ is regular and $e$ closed, there is a closed
neighbourhood $V$ of $y$ with $e^{-1}(V)\subs\bigcup\{ U_z\: |\: z\in e^{-1}(y)
\}$. Hence
$$U_x\cap e^{-1}(V)=e^{-1}(V)\setm\bigcup\bigl\{ U_z\: |\: z\in
e^{-1}(y)\setm\{ x\}\bigr\}$$
is a closed neighbourhood of $x$. Thus $C$ is regular.

\vspc
Let $F\subs C$ be finite. Then $\overline{C(e(F))}$ is compact. By
Proposition~\ref{p10}, the fibers of $e$ are compact. Hence
$e^{-1}(\overline{C(e(F))})$ is compact by \cite{Bou}, chap. I.10,
Proposition~6. Furthermore, $e^{-1}(\overline{C(e(F))})$ is convex with respect
to $\VVV_e$. Therefore, the closed subset $\overline{C(F)}$ of
$e^{-1}(\overline{C(e(F))})$ is compact. This proves (G1) for $C$.

\vspc
Next let $e'\colon C'\ra C$ be an \'etale map with $C'$ compact which is
generated by $\{ x,y\}\subs C'$. Then $ee'$ is \'etale and generated by
$\{ x,y\}$. Hence $C'$ is minimal among the connected sets $B\subs C'$ with
$x,y\in B$. Thus $C$ satisfies (G2).

\vspc
Finally, let $S(x,E):=\bigcup\{ C(x,z)\:|\:z\in E\}$ be a closed star in some
$U\in\VVV_e$. Since $C$ is regular, we find a closed convex neighbourhood
$U'\subs U$ of $x$. By Proposition~\ref{p3}, this implies that $S(x,E)\cap U'$
is a star in $U$ which is closed in $C$.
Therefore, $e(S(x,E)\cap U')$ is a closed star in some $V\in\VVV$. So $E$ is
finite, which proves that $C$ is a geodesic $q$-manifold. \bx

Now we are ready to prove our main result which essentially states that the
image of an \'etale map is weakly convex in the following sense
(cf. \cite{BOR1}, Definition~2.16).
\begin{Definition}\pt\label{d9}
\rm Let $Y$ be a geodesic manifold. We call a subset $A\subs Y$ {\em weakly
convex} if every pair of points $x,y\in A$ can be connected by a geodesic.
\end{Definition}

The following theorem extends previous versions of the Lokal-global-Prinzip
for convexity of maps (see \cite{CDM, HNP, BOR1,BOR2}).

\begin{thm}\pt\label{t3}
Let $f\colon X\ra Y$ be a locally convex continuous map from a connected
topological
space $X$ to a geodesic $q$-manifold $Y$. Assume that $f^{\#}$ is closed.
Then $f(X)$ is weakly convex.
\end{thm}

\pf Let $\VVV$ be the atlas of $Y$. By Proposition~\ref{p7}, the map $f^{\#}$
again locally convex, and Proposition~\ref{p5} implies that $f^{\#}$ is
\'etale. By Theorem~\ref{t2}, it follows that $X^f$ is a geodesic manifold.
For $z,z'\in X^f$, Proposition~\ref{p8} shows that there is a connecting
simple line path $L$ between $z$ and $z'$. Theorem~\ref{t1} shows that $L$ can
be straightened to a convex simple geodesic $C$. Thus $f^{\#}|_C\colon C\ra Y$
is a geodesic between $f^{\#}(z)$ and $f^{\#}(z')$. Hence $f(X)$ is
weakly convex. \bx

In the special case where $f$ is an inclusion $X\hra Y$, the preceding proof
yields

\vspace{3mm}
\noindent {\bf Corollary.} \it Let $C$ be a closed connected locally convex
subset of a geodesic manifold $Y$. Then $C$ is weakly convex.       \rm

\vspace{3mm}
\pf By Example~3, $C$ is a geodesic manifold, and $C\hra Y$ is \'etale. As in
the proof of Theorem~\ref{t3}, this implies that $C$ is weakly convex. \bx

\noindent {\bf Remarks. 1.} If $f$ is closed, then $f^{\#}$ is closed. However,
the latter condition is much weaker. For example, if $f$ is a local
homeomorphism, then $f^{\#}$ is identical, but $f$ need not be closed.

{\bf 2.} The preceding corollary extends Klee's generalization
of a classical result due to Tietze \cite{Ti} and Nakajima (Matsumura)
\cite{Na}. Klee's theorem \cite{Klee} states that the above corollary holds in
a locally convex topological vector space $Y$. Note that the usual proof of
Klee's theorem rests on the linear structure of $Y$, while the corollary of
Theorem~\ref{t3} merely depends on a local convexity structure in the sense of
Definition~\ref{d4}.

\end{document}